\documentstyle[11pt]{article}
\date{}

\def\l{{\cal L}_{{\cal B}}(X)}

\def\s1{{\cal S}_1 (X)}

\def\f{\varphi}

\def\ben{\begin{enumerate}}
\def\een{\end{enumerate}}

\def\zR{{\rm I}\!{\rm R}}

\def\zZ{{\rm Z}\!\!{\rm Z}}
\def\zC{{\rm C}\!\!\!\vrule height 6.4pt depth -0.4pt width
1pt \ }

\def\m{{\cal M}}
\def\1{\vec{1}}
\def\l2{L^2(\m,\tau)}

\def\hrarr^#1_#2{ \mathrel{
\mathop{\hbox to .3in{\rightarrowfill}}
\limits^{\scriptstyle#1}_{\scriptstyle#2}  }}

\def\hlarr^#1_#2{ \mathrel{
\mathop{\hbox to .5in{\leftarrowfill}}
\limits^{\scriptstyle#1}_{\scriptstyle#2}  }}

\newtheorem{lemm}{Lemma}[section]

\newtheorem{prop}[lemm]{Proposition}
\newtheorem{coro}[lemm]{Corollary}
\newtheorem{rema}[lemm]{Remark}

\newtheorem{exam}[lemm]{Example}

\newenvironment{proof}[1]{
  \trivlist \item[\hskip \labelsep{\bf #1}]}{\hfill\mbox{$\Box$}
  \endtrivlist}

\begin{document}

\title{ Connes' metric for states in group algebras}

\author{ Esteban Andruchow and Gabriel Larotonda}
\maketitle




\vskip0.5cm

\abstract{ In this article we follow the main idea of  A. Connes for the construction of a metric in the state space of a $C^*$-algebra. We focus in the reduced algebra of a discrete group $\Gamma$, and prove some equivalences and relations between two central objects of this category: the word-length growth (connected with the degree of the extension of $\Gamma$ when the group is an extension of $\zZ$ by a finite group), and the topological equivalence between the $\omega^*$ topology and the one introduced with this metric in the state space of $C_r^*(\Gamma)$.}

\bigskip


\bigskip
\section{Introduction}
In \cite{con1} and \cite{con2}, A. Connes introduced what he called non commutative metric spaces, which consist of a triples $({\cal A},D,H)$ where ${\cal A}$ is a C$^*$-algebra, acting on the Hilbert space $H$, and $D$ is an unbounded operator in H, called the Dirac operator, satisfying 
\begin{itemize}
\item $(D^2+1)^{-1}$ is compact
\item the set $\{a\in{\cal A}:[D,a]\;is\; bounded\;\}$ is norm-dense in ${\cal A}$
\end{itemize}

We are interested in the case when $\Gamma$ is a discrete group with identity element $e$ and the algebra ${\cal A}$ is the reduced C$^*$-algebra $C^*_r(\Gamma)$. The Hilbert space is $\ell^2(\Gamma)$, with $C^*_r(\Gamma)$ acting as left convoluters (i.e. the left regular representation). The Dirac operator is defined in terms of a length function on $\Gamma$. A length function is a map $L:\Gamma \to \zR_+$ satisfying
\begin{enumerate}
\item
$L(gh)\le L(g)+L(h)$ for all $g,h \in \Gamma$.
\item
$L(g^{-1})=L(g)$ for all $g\in \Gamma$.
\item
$L(e)=0$.
\end{enumerate}
If $\Gamma$ is given by generators and relations, the prototypical length function is the map which assigns to each word its (minimal) length. We shall fix this data $L$, and we will make the further assumption that the sets
$$
\{g\in \Gamma: L(g)\le c\}
$$
are finite for any $c>0$. The Dirac operator \cite{con2} is then defined as follows:
$$
D(\delta_g)=L(g)\delta_g
$$
where $\{ \delta_g: g\in \Gamma\}$ is the canonical orthonormal basis of $\ell^2(\Gamma)$. As is custom, we shall denote by $\lambda_g$ the element $\delta_g$ regarded as an  operator in $\ell^2(\Gamma)$. The metric (of the non commutative metric space) is defined in the state space ${\cal S}(C^*_r(\Gamma))$ of $C^*_r(\Gamma)$ by means of the formula
$$
d(\psi,\f)=\sup\{|\psi(a)-\f(a)| : a\in C^*_r(\Gamma) \hbox{ with } \|[a,D]\|\le 1\}.
$$
Here $[\ , \ ]$ denotes the usual conmutator of operators. This $d$ is not necessarily finite. In this note we study situations in which it is finite, and consider a problem posed by M. Rieffel, asking under which assumptions the metric thus defined induces on the state space a topology which is equivalent to the $w^*$ topology.

The basic example of this situation, which
even justifies the name "non commutative metric space", occurs when ${\cal A}$ is $C(M)$, the algebra of continuous functions on a spin manifold $M$ \cite{con2}, \cite{tesis}. M Rieffel found \cite{rifel} a natural triple associated to the noncommutative tori. Also he pointed out that one can find a positive answer for matrix algebras.

In this note we consider this problem for group algebras arising from discrete groups and  triples arising from length functions. Instead of dealing with the $d$ metric directly, we refer it to two metrics, $d_\infty$ and $d_2$, related with the asymptotic behaviour of the family $\{\frac{1}{L(g)}: e\neq g\in \Gamma\}$:
\begin{equation} \label{dinfinito}
d_\infty(\f,\psi)=\sup_{g\neq e}\frac{|\f(\lambda_g)-\psi(\lambda_g)|}{L(g)} ,
\end{equation}
and
\begin{equation}\label{deledos}
d_2(\f,\psi)=\bigl(\sum_{g\neq e}\frac{|\f(\lambda_g)-\psi(\lambda_g)|^2}{L(g)^2} \bigr)^{1/2}.
\end{equation}
First note that $d_\infty$ is a well defined metric and that $d_\infty(\f,\psi)\le d(\f,\psi)$. The first fact is apparent. To prove the second, note that $\|[D,\lambda_g]\|=L(g)$, and therefore
$$
d_\infty(\f,\psi)=\sup_{a=\frac{1}{L(g)}, g\neq e}|\f(a)-\psi(a)|\le d(\f,\psi) .
$$
Also note that $d_2$ may fail to be finite. Indeed, consider
$\Gamma=\zZ\times \zZ$. Then the family $\{\frac{1}{L(g)}: g\neq
e\}$ does not belong to $\ell^2(\zZ\times \zZ)$.  Consider the
positive definite functions $f(g)=1$ for all $g$ and
$h=\delta_e$. These functions induce states $\f_f$ and $\f_h$ on
$C^*_r(\zZ\times \zZ)$ satisfying $\f_f(\lambda_g)=f(g)$ and
$\f_h(\lambda_g)=h(g)$. It follows that
$$
d_2(\f_f,\f_h)=\sum_{e\neq g\in \zZ\times
\zZ}\frac{1}{L(g)^2}=\infty .
$$
Denote by $K(\Gamma)$ the group algebra of $\Gamma$, i.e. the set
of elements of the form $\sum_{g\in F}\alpha_g \lambda_g$, where
$\alpha_g\in \zC$ and $F\subset \Gamma$ is a finite set.

\begin{lemm}
$d_\infty$ is a metric in ${\cal S}(C^*_r(\Gamma))$ which induces
a topology equivalent to the $w^*$-topology.
\end{lemm}
\begin{proof}{Proof.}
If $d_\infty(\f_n,\f)\to 0$, then clearly $\f_n(\lambda_g)\to
\f(\lambda_g)$ for all $g\neq e$. Since $\f_n$,$\f$ are states,
$\f_n(\lambda_e)=\f_n(1)=1=\f(\lambda_e)$. It follows that
$\f_n(a)\to \f(a)$ for all $a\in K(\Gamma)$. Since $\f_n$, $\f$
have their norms bounded (by $1$), and since $K(\Gamma)$ is dense
in $C^*_r(\Gamma)$, it follows that $\f_n\to \f$ in the $w^*$
topology. Conversely, suppose that $\f_n(a)\to \f(a)$ for all
$a\in C^*_r(\Gamma)$ and fix $\epsilon>0$. Let
$F=\{g\in\Gamma: L(g)<4/\epsilon\}$, which is a finite
set, say $F=\{g_1,...,g_k\}$. If $g\in F$, one
has
$$
\frac{|\f_n(\lambda_g)-\f(\lambda_g)|}{L(g)}\le
\frac{|\f_n(\lambda_g)|+|\f(\lambda_g)|}{4/\epsilon}\le\epsilon/2
.
$$
On the other hand, there exists $n_0$ such that for all $n\ge
n_0$, 
\newline
$|\f_n(\lambda_{g_i})-\f(\lambda_{g_i})|<\frac{\epsilon}{2}
\min\{L(g_1),...,L(g_k)\}$, for $i=1,...,k$. It follows that
$$
\frac{|\f_n(\lambda_{g_i})-\f(\lambda_{g_i})|}{L(g_i)}<\epsilon/2
.
$$
Therefore, if $n\ge n_0$, then
$$
\sup_{g\neq
e}\frac{|\f_n(\lambda_g)-\f(\lambda_g)|}{L(g)}=d_\infty(\f_n,\f)\to
0 .
$$
\end{proof}
\section{Comparison between $d$,$d_\infty$ and $d_2$}

Here we establish the basic inequality for these metrics, namely
$d_\infty\le d \le d_2$.

\begin{lemm}
Let $a= \sum_{g\in F}\alpha_g \lambda_g \in K(\Gamma)$, then
$$
\|[D,a]\|\ge \bigl(\sum_{g\in F} |\alpha_g|^2 L(g)^2\bigr)^{1/2} .
$$
\end{lemm}
\begin{proof}{Proof.}
Note that
$$
[D,a]\delta_e=\sum_{g\in F} \alpha_g
[D,\lambda_g]\delta_e=-\sum_{g\in F}\alpha_g L(g) \delta_g ,
$$
because $D\lambda_g\delta_e=D\delta_g=L(g)\delta_g$, and in
particular $D\delta_e=0$. Therefore
$$
\|[D,a]\delta_e\|_2^2=\sum_{g\in F} |\alpha_g|^2 L(g)^2 .
$$
It follows that $\|[D,a]\|\ge
\|[D,a]\delta_e\|_2=\bigl(\sum_{g\in F} |\alpha_g|^2
L(g)^2\bigr)^{1/2}$.
\end{proof}
\begin{prop}
$$
d_\infty(\f,\psi)\le d(\f,\psi)\le d_2(\f,\psi) .
$$
\end{prop}
\begin{proof}{Proof.}
Pick $a=\sum_{g\in F}\alpha_g \lambda_g\in K(\Gamma)$, with
$\|[D,a]\|\le 1$ (note that for any $a\in K(\Gamma)$, $[D,a]$ is
a bounded operator). Then
$$
|\f(a)-\psi(a)|=|\sum_{g\in
F}\alpha_g(\f(\lambda_g)-\psi(\lambda_g))|=|\sum_{e\neq g\in F}
\alpha_g L(g)\frac{(\f(\lambda_g)-\psi(\lambda_g))}{L(g)}| ,
$$
which by the Cauchy-Schwartz inequality is less than or equal to
$$
\bigl(\sum_{e\neq g\in
F}|\alpha_g|^2L(g)^2\bigr)^{1/2}\bigl(\sum_{e\neq g\in
F}\frac{|\f(\lambda_g)-\psi(\lambda_g)|^2}{L(g)^2}\bigr)^{1/2}\le
\|[D,a]\|d_2(\f,\psi)\le d_2(\f,\psi).$$ 

The proof finishes by
observing that the set of elements $a\in K(\Gamma)$ with
$\|[D,a]\|\le 1$ is dense among  elements  $b\in C^*_r(\Gamma)$ with $\|[D,b]\|\le 1$.
Indeed, let $b=\sum_{g\in \Gamma}\beta_g \lambda_g \in
C^*_r(\Gamma)$ with $\|[D,b]\|\le 1$. For  finite sets $F\subset
\Gamma$, the truncated elements $b_F=\sum_{g\in F}\beta_g\lambda_g
\in K(\Gamma)$ converge in norm to $b$. Clearly also the
(bounded) commutants $[D,b_F]$ converge in norm to $[D,b]$.

\noindent Denote by $N_F=\|[D,b]\| \|[D,b_F]\|^{-1}$ (after droping the elements
$b_F$ with $[D,b_F]=0$). Then $N_F b_F$ lies in $K(\Gamma)$, the commutants $[D,N_F b_F]$ have norm less than or equal to $1$,
and converge to $b$.
\end{proof}
We emphasize that $d_2$ might be infinite. It would be finite if for example the family $\{\frac{1}{L(g)}: e\neq g\in \Gamma\}$ would lie in $\ell^2(\Gamma)$. This imposes a strong condition on $\Gamma$, namely, that the group $\Gamma$ has linear growth (polinomial growth with degree $1$), see \cite{gromov} and \cite{con2}.  
This means, that there exists constants $k,l$ such that $\#\{g\in \Gamma: L(g)\le c\}\sim kc+l$. 
\begin{exam}
Let us consider the following examples, of groups $\Gamma$ wich satisfy that the family $\{\frac{1}{L(g)}: e\neq g\in \Gamma\}$ lies in $\ell^2$.
\begin{enumerate}
\item
Let $\Gamma=\zZ$. Here the length function is $L(m)=|m|$, $m\in \zZ$. The group C$^*$-algebra equals in this case $C(S^1)$.
\item 
Let $\Gamma$ be a finite extension of $\zZ$, i.e. a group $\Gamma$ which has a copy of $\zZ$ inside, as a normal subgroup, and the quotient ${\cal F}=\Gamma / \zZ$ is finite. Then, as a set, $\Gamma$ is $\zZ\times {\cal F}$. Let ${\cal F}=\{f_1,...,f_n\}$. Then the classes of $(1,f_1),..., (1,f_n)$ (i.e. these elements regarded as elements of $\Gamma$) are generators for $\Gamma$. Let us consider the length function $L$ given by word length with respect to this set of generators. Note that for this $L$, there are at most $2n$ elements of $\Gamma$ with any given length. It follows that  $\{\frac{1}{L(g)}: e\neq g\in \Gamma\}$ lies in $\ell^2$. The (reduced) $C^*$-algebra of such $\Gamma$ can be computed. They consist of algebras of $n\times n$ matrices with entries in $C(S^1)$, see chapter VIII of $\cite{dav}$ for a complete descrition of this computation. Let us point out  two special cases of this type
\begin{enumerate}
\item
$\Gamma=\zZ\times {\cal F}$ with the usual product for pairs. In this case the $C^*$-algebra is $C^*_r(\zZ\times {\cal F})\simeq C(S^1)\otimes C^*_r({\cal F})$.  The algebra $C^*_r({\cal F})$ is finite dimensional, therefore in this case $C^*_r(\Gamma)$ consists of a direct sum of full matrix algebras with entries in $C(S^1)$. In particular, if ${\cal F}=S_k$ the group of permutations of order $k$, then $C^*_r(\Gamma)=M_k(C(S^1))$.
\item
Consider the (unique) nontrivial automorphism of $\zZ$, $\theta(m)=-m$. Then one has a $\zZ_2$ extension of $\zZ$, $\Gamma=\zZ \times_\theta \zZ_2$, and the corresponding C$^*$-algebra  $C^*_r(\Gamma)$ is the cross product $C(S^1)\times_\theta \zZ_2$, which identifies with the algebra of $2\times 2$ matrices with entries in $C(S^1)$ of the form
$$
\left(
\begin{array}{lr}

f(z) & g(z) \\ f(\overline{z}) & g(\overline{z}) 
 
\end{array}
\right),
$$
where $f$ and $g$ are continuous functions in $S^1$.

\end{enumerate}
\end{enumerate}
\end{exam}

\vspace*{1cm}

\begin{prop} 
If $\Gamma$ has a length function $L$ which satisfies that $\{\frac{1}{L(g)}: e\neq g\in \Gamma\}$ is square summable, then the metric $d$ is well defined (is finite) and induces on the state space of $C^*_r(\Gamma)$ the $w^*$ topology.
\end{prop}
\begin{proof}{Proof.}
Since $\{\frac{1}{L(g)}\}\in\ell^2$, $d\le d_2<\infty$. By the above results it suffices to prove that if a sequence $\f_n$ converges to $\f$ in the $w^*$ topology, then it converges in the $d$ metric. We claim that it converges in the $d_2$ metric. Fix $\epsilon >0$. There exists a finite set $F=\{g_1,...,g_k\}$ such that 
$$
\bigl( \sum_{g\in \Gamma - F} \frac{|\f_n(\lambda_g)-\f(\lambda_g)|^2}{L(g)^2}\bigr)^{1/2}
\le 2\bigl( \sum_{g\in\Gamma -F} \frac{1}{L(g)^2}\bigr)^{1/2} <\epsilon/2 .
$$
Put $c=(\sum_{i=1}^k \frac{1}{L(g_i)^2})^{1/2}$. There exists $n_0$ such that for all $n\ge n_0$, one has $|\f_n(\lambda_{g_i})-\f(\lambda_{g_i})|<\epsilon /2c$. Therefore
$$
\bigl( \sum_{i=1}^k \frac{|\f_n(\lambda_{g_i})-\f(\lambda_{g_i})|^2}{L(g_i)^2}\bigr)^{1/2} <\epsilon/2 .
$$
Then $d_2(\f_n,\f)\to 0$.
\end{proof}

\bigskip

\begin{coro}

If $\Gamma$ is a finite extension of $\zZ$, then in $S\left(C_r^*(\Gamma)\right)$ the $d$ metric is well defined and induces the $\omega^*$-topology.

\end{coro}

\bigskip

\section{Normal states which are bounded with respect to the trace}

We shall prove that the metric $d$ is finite on the set of normal states of $C^*_r(\Gamma)$ which are bounded with respect to the trace of $C^*_r(\Gamma)$, i.e. the states $\f$ which extend to normal states of the Von Neumann algebra ${\cal L}_\Gamma$ of $\Gamma$, and verify that there exists a constant $\kappa>0$ such that
$$
\f(a^*a)\le \kappa \tau(a^*a),
$$
or shortly, $\f \le \kappa \tau$.
Recall that the trace $\tau$ is given by $\tau(a)=\langle a\delta_e,\delta_e\rangle$. There is a Radon-Nykodim derivative for all such $\f$ \cite{araki}. Namely, there exists an element $\rho_\f \ge 0$ in ${\cal L}_\Gamma$ such that
$$
\f(a)=\tau(\rho_\f a) \ , \ \ \hbox{ with } \|\rho_\f\|\le \kappa^{1/2} .
$$
Denote by ${\cal S}_\kappa$ the set
$$
{\cal S}_\kappa=\{\f\in {\cal S}({\cal L}_\Gamma): \f\le \kappa \tau \} .
$$
First note that a state which lies in ${\cal S}_\kappa$ is necessarily normal. Indeed, let $\{p_i:i\in I\}$ be an arbitrary family of pairwise orthogonal projections in ${\cal L}_\Gamma$. Fix $\epsilon>0$ and let $J\subset I$ be a finite set such that $\tau(\sum_{i\in I-J}p_i)=\sum_{i\in I-J}\tau(p_i)<\epsilon/\kappa$. Then $\f(\sum_{i\in I-J}p_i)<\epsilon$. Therefore $0\le \f(\sum_{i\in I}p_i)=\sum_{j\in J}\f(p_j)+\f(\sum_{i\in I-J}p_i)\le \sum_{j\in J}\f(p_j)+\epsilon$. That is, $\sum_{i\in I}\f(p_i)=\f(\sum_{i\in I}p_i)$, and $\f$ is normal.

Also it is apparent that ${\cal S}_\kappa$ is $w^*$ compact and convex.

\begin{prop} 
The metrics $d$ and $d_2$ are well defined on ${\cal S}_\kappa$ and induce  the $w^*$ topology.
\end{prop}
\begin{proof}{Proof.}
Note that if $\f\in {\cal S}_\kappa$ then 
$$
\f(\lambda_g)=\tau(\rho_\f \lambda_g)=\langle \rho_\f \delta_g,\delta_e\rangle=\rho_\f(g^{-1}),
$$
where $\rho_\f(g^{-1})$ denotes the $g^{-1}$-coordinate of $\rho_\f$ regarded as an element of $\ell^2(\Gamma)$. In particular, it follows that the family $\{\f(\lambda_g): g\in \Gamma\}$ is square summable. Moreover, 
$$
\bigl(\sum_{g\in \Gamma} |\f(\lambda_g)|^2\bigr)^{1/2} =\|\rho_\f\|_2 \le \|\rho_\f\|\le \kappa^{1/2}.
$$
It follows that if $\f,\psi\in {\cal S}_\kappa$, them 
$$
d(\f,\psi)\le d_2(\f,\psi)=\bigl(\sum_{e\neq g\in \Gamma} \frac{|\f(\lambda_g)-\psi(\lambda_g)|^2}{L(g)^2}\bigr)^{1/2}\le 2\kappa .
$$
Suppose now that $\f_n\to \f$ in the $w^*$ topology. Fix $\epsilon>0$. Then the  set $F=\{g\in \Gamma: L(g)\le 2\kappa/\epsilon\}$ is finite. Say $F=\{g_1,...,g_n\}$. If $g$ lies outside $F$ one has
$$\bigl(\sum_{e\neq g\in \Gamma-F} \frac{|\f(\lambda_g)-\psi(\lambda_g)|^2}{L(g)^2}\bigr)^{1/2}\le\frac{\epsilon}{2\kappa} \bigl(\sum_{g\in \Gamma -F}|\f_n(\lambda_g)-\f(\lambda_g)|^2\bigr)^{1/2}\le \epsilon .
$$
Let $C=\sum_{i=1}^n \frac{1}{L(g_i)^2}$. There exists $n_0$ such that if $n\ge n_0$ then
$$
|\f_n(\lambda_{g_i})-\f(\lambda_{g_i})|<\epsilon/C \ \hbox{ for } i=1,...,n.
$$
It follows that $d(\f_n,\f)\le d_2(\f_n,\f)<\epsilon$ if $n\ge n_0$.
\end{proof}
\begin{rema}
\begin{enumerate}
\item
The first part of the proof in fact shows that if $\f$ and $\psi$ are normal states of ${\cal L}_\Gamma$ whose Radon-Nykodim derivatives with respect to the trace $\tau$ lie in $\ell^2(\Gamma)$, then $d(\f,\psi)\le d_2(\f,\psi)<\infty$ .
\item
It is apparent that the metrics $d$ and $d_2$ are also finite on the set
$\cup_{\kappa>0} {\cal S}_\kappa$
\end{enumerate}
\end{rema}

\noindent
\begin{tiny}
\noindent
\hspace*{-.05cm}Instituto de Ciencias\\
Universidad Nacional de Gral. Sarmiento\\
J. M. Gutierrez 1150\\
(1613) Los Polvorines\\
Argentina \\
e-mail: eandruch@ungs.edu.ar, glaroton@ungs.edu.ar
\end{tiny}

\vskip2cm

{
}

\end{document}